\documentclass[letterpaper, 10 pt, conference]{ieeeconf} 
\makeatletter
\let\NAT@parse\undefined
\makeatother

\makeatletter
\renewcommand{\paragraph}[1]{%
  \par\vspace{0.4em}%
  \noindent\textbf{#1.}%
}
\makeatother

\let\labelindent\relax

\title{\LARGE \bf Accelerating Full-Scale Nonlinear Model Predictive Control via Surrogate Dynamics Optimization}
\usepackage[utf8]{inputenc}
\usepackage[T1]{fontenc}      
\usepackage{cite}
\usepackage{amsmath}               
\allowdisplaybreaks
\usepackage{amssymb}               
\usepackage{amsfonts}              
\usepackage{amsbsy}
\usepackage{bbm}
\usepackage{tikz}                  
\usepackage{multicol}
\usepackage{graphicx}
\usepackage{threeparttable}
\usepackage{soul,derivative,stmaryrd,mathtools,enumerate,enumitem,mathtools}
\usepackage{comment}
\usepackage{algorithm}
\usepackage{algorithmicx}
\usepackage{algpseudocode}
\usepackage{xfrac}
\usepackage{multirow}
\usepackage{booktabs}

\usetikzlibrary{arrows,positioning,shapes}
\usetikzlibrary{calc}
\usetikzlibrary{shapes.multipart}
\usetikzlibrary{arrows.meta}
\usetikzlibrary{decorations.pathreplacing}
\usetikzlibrary{arrows,intersections}
\usetikzlibrary{shapes,arrows}

\newcommand{\ol}[1]{#1}

\newcommand{\mbb}[1]{\mathbb{#1}}
\newcommand{\mbf}[1]{\mathbf{#1}}
\newcommand{\mcal}[1]{\mathcal{#1}}
\newcommand{\mrm}[1]{\mathrm{#1}}

\newcommand{\R}{\mbb{R}}

\renewcommand{\P}{\mbb{P}}
\newcommand{\E}{\mbb{E}}

\newcommand{\N}{\mbb{N}}

\newcommand{\U}{\mcal{U}}
\newcommand{\AO}{\mrm{AO}}
\newcommand{\I}{\mrm{I}}
\newcommand{\Pn}{\mrm{P}}

\newcommand{\Xe}{\mrm{X}}
\newcommand{\Tin}{\mrm{T}_\mrm{in}}
\newcommand{\Tref}{\mrm{T}_\mrm{ref}}
\newcommand{\Cb}{\mrm{C}_\mrm{b}}
\newcommand{\h}{\mrm{h}_{\mrm{cr}}}
\newcommand{\nz}{n_{\mrm{z}}}

\DeclareMathOperator*{\argmin}{arg\,min}
\newcommand{\sss}{\scriptscriptstyle}

\DeclareMathSymbol{\shortminus}{\mathbin}{AMSa}{"39}
\renewcommand{\between}[2]{%
  \mathord{\text{\scriptsize[$#1\!\!:\!\!#2$]}}%
}

\newtheorem{remark}{Remark}
\newtheorem{proposition}{Proposition}
\usepackage{subcaption} 

\definecolor{couleur1}{HTML}{0f4c5c}
\definecolor{couleur2}{HTML}{fb8b24}
\definecolor{couleur3}{HTML}{5f0f40}
\usepackage[colorlinks=true,citecolor=blue]{hyperref}

\IEEEoverridecommandlockouts                              

\author{
Perceval Beja-Battais$^{1,2}$
 \thanks{$^{1}$Université Paris Saclay, ENS Paris Saclay, CNRS, Centre Borelli, 91190, Gif-sur-Yvette, France 
        {\tt\small name.surname@ens-paris-saclay.fr}}
\quad
 Guillaume Dupré$^{2}$ \quad
 \thanks{$^{2}$Framatome, Reactor Process Department, 92400, Courbevoie, France 
        {\tt\small name.surname@framatome.com} \newline
        This work has been submitted to the IEEE for possible publication. Copyright may be transferred without notice, after which this version may no longer be accessible.}
Alain Grossetête$^{2}$
\quad
Nicolas Vayatis$^{1}$
}

\begin{document}

\maketitle
\thispagestyle{empty}
\pagestyle{empty}
\begin{abstract}                
Driven by advances in hardware and software technologies, nonlinear model predictive control (NMPC) has gained increasing adoption in both industry and academia over the past decades. 
However, its practical deployment is often limited by the computational cost of simulating the embedded process model, especially for high-dimensional, multi-time-scale, or nonlinear systems commonly found in real-world applications.
Thus, this paper introduces Surrogate Dynamics Optimization (SDO), a warm-start framework for full-scale NMPC to address the limitation of standard initialization strategies.
The approach relies on a machine learning surrogate model to solve a lightweight auxiliary problem that approximates the original one.
The methodology is reproducible and compatible with in-house simulation and optimization tools, a key consideration in industrial contexts.
Data efficiency of SDO, as well as the impact of surrogate design on the overall performance, are evaluated through a non-trivial simulation case study: 24-hour optimal load-following control of a pressurized water reactor. 
The results show consistent improvements in NMPC convergence speed within a fixed computational budget, while reducing training data generation costs by two orders of magnitude compared to behavior cloning.
\end{abstract}



\section{Introduction}\label{sec:Introduction}

Model Predictive Control (MPC) is a widely used control methodology that has been successfully applied in many industrial sectors~\cite{qin2003survey}. 
It relies on the online, repeated solution of an Optimal Control Problem (OCP), which minimizes a cost function subject to model dynamics and additional constraints. 
Nonlinear MPC (NMPC) extends this approach to OCPs with nonlinear dynamics, cost functions, or constraints. 
Despite its higher implementation complexity~\cite{rao2009survey, diehl2009efficient}, NMPC has attracted growing interest in industry~\cite{schwenzer2021review}, as most systems are inherently nonlinear. 
However, several obstacles continue to limit its practical deployment~\cite{fiedler2023mpc}, despite recent advances in real-time NMPC schemes~\cite{wolf2016fast}.
A common bottleneck is the computational time required to solve the online OCP, particularly for long-horizon problems.

To mitigate this bottleneck, warm-start strategies can be implemented to provide good initial guesses for the decision variables.
Traditional approaches, such as shift initialization~\cite{diehl2009efficient}, exploit the receding-horizon structure of the OCP.
However, as noted by~\cite{li2024faster}, these simple strategies become ineffective when the cost function, constraints, or operating conditions change abruptly, since the previous optimal trajectory may no longer be valid.

With the growing availability of data and computational resources, Machine Learning (ML) has emerged as another effective means of reducing online computational costs.
Many recent applications have been reported in the MPC literature~\cite{hewing2020learning}, typically to learn either the system dynamics~\cite{hossain2023machine} or the control policy~\cite{vaupel2020accelerating, li2024faster}, both of which can be used to provide high-quality initial guesses.
This naturally motivates the use of learning-based approaches to provide high-quality initializations for long-horizon NMPC. 
Still, such methods also tend to exhibit critical limitations in the setting of this work, which are analyzed in more detail in Sec.~\ref{sec:pb_statement} once the problem structure is established.
Furthermore, to the best of our knowledge, existing learning-based warm-start strategies are largely problem-specific and tailored to explicit applications, such as autonomous driving~\cite{bouzidi2024learning} or legged robotics~\cite{kim2025real}.
To address these challenges, this paper investigates how ML can accelerate NMPC in settings where full-scale model evaluations constitute the primary computational bottleneck, while preserving the structure of the classical control pipeline.
The main contributions and results of this paper are as follows:

\begin{figure}
    \centering
    \includegraphics[width=\linewidth]{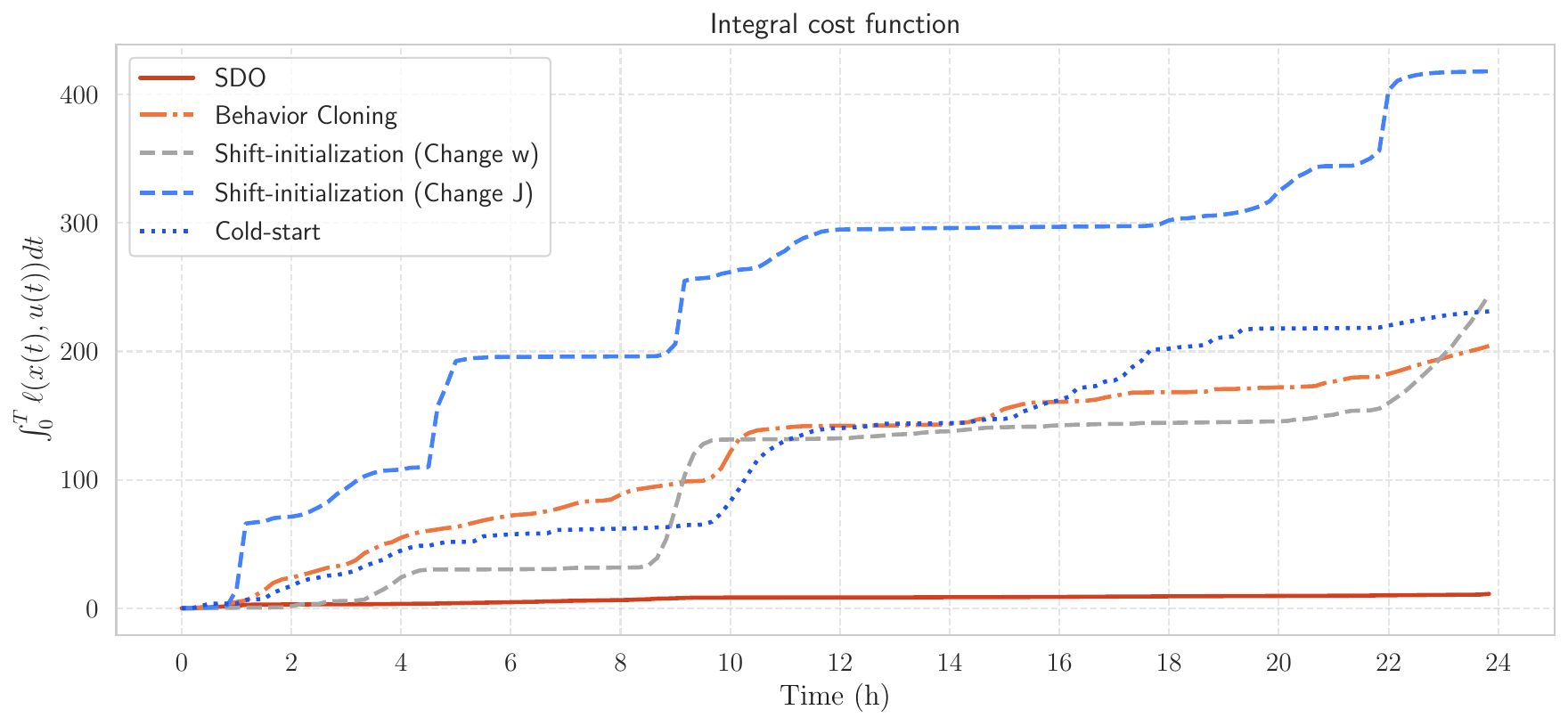}
    \caption{Comparison of the objective cost achieved after a fixed computational budget for different warm-start strategies on the industrial case study, showing faster convergence of SDO (see Sec.~\ref{sec:experiments}).}
    \label{fig:main_result}
\end{figure}



\paragraph{A novel methodology for long-horizon NMPC with expensive model evaluations}
Surrogate Dynamics Optimization (SDO) is specifically designed for single-shooting transcriptions and high-dimensional control spaces (e.g., $\geq 100$). 
In Section~\ref{sec:methodology}, we introduce a framework that combines data-driven modeling of the system dynamics with an optimization routine to produce structured initial control sequences. 

\paragraph{An offline pipeline with reduced sampling requirements} 
The surrogate model is trained to approximate the system dynamics rather than optimal control trajectories.
Thus, the required dataset can be generated significantly faster than in imitation learning-based approaches. 
This makes the method suitable for settings where full-scale optimal trajectory generation is expensive.
Moreover, we experimentally demonstrate the effectiveness of the framework when the surrogate model is trained on fewer data points.

\paragraph{Safe integration within the full-scale NMPC pipeline and validation on an industrial benchmark}
The surrogate-based solution is systematically refined through re-optimization in the original NMPC formulation, ensuring feasibility with respect to the full-scale model and constraints. 
The approach is validated on a realistic industrial case study involving 24-hour load-following control of Pressurized Water Reactors (PWRs). 
Under a fixed computational budget, SDO significantly improves the convergence rate of the NMPC solver and mitigates worst-case performance degradation compared to shift initialization and behavior cloning strategies.




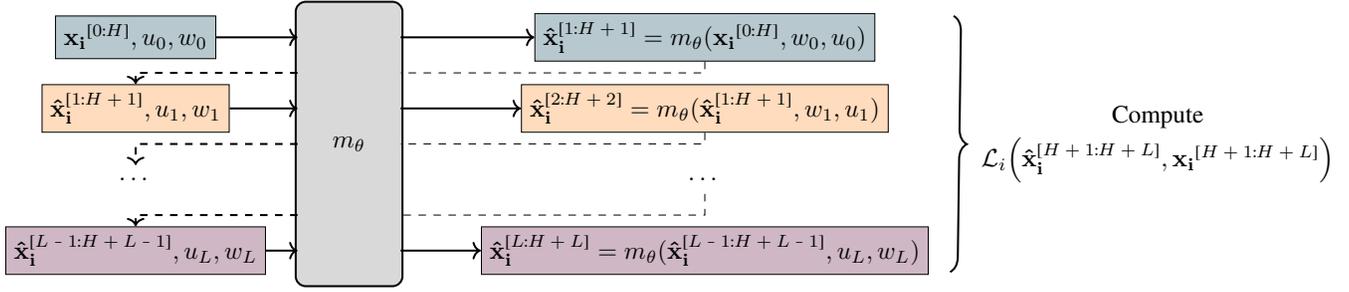
\begin{figure*}
    \centering
    \resizebox{1\linewidth}{!}{\begin{tikzpicture}[
    xnode/.style={draw, rectangle, minimum width=1cm, minimum height=0.6cm},
    pred/.style={draw, rectangle, dashed,minimum width=1cm, minimum height=0.6cm, node distance=3cm},
    arrow/.style={->, thick},]
    
    \node[xnode, fill=couleur1!30] (x0H) at (0,0) {$\mbf{x_i}^{{\tiny\between{0}{H}}}, u_0, w_0$};
    \node[xnode, fill=couleur2!30] (x1Hp1) at (0,-1) {$\mbf{\hat{x}_i}^{{\tiny\between{1}{H+1}}}, u_1, w_1$};
    \node (dots1) at (0, -2) {$\cdots$};
    \node[xnode, fill=couleur3!30] (xLHpL) at (0,-3) {$\mbf{\hat{x}_i}^{{\tiny\between{L\shortminus1}{H+L\shortminus1}}}, u_L, w_L$};
    \node[draw, thick, rounded corners, align=center, minimum width=1.5cm, minimum height=4cm, fill=gray!30] (model) at (3, -1.5) {$m_\theta$};

    \node[xnode, fill=couleur1!30] (out0H) at (8,0) {$\mbf{\hat{x}_i}^{{\tiny\between{1}{H+1}}} = m_\theta(\mbf{x_i}^{{\tiny\between{0}{H}}}, w_0, u_0)$};
    \node[xnode, fill=couleur2!30] (out1Hp1) at (8,-1) {$\mbf{\hat{x}_i}^{{\tiny\between{2}{H+2}}} = m_\theta(\mbf{\hat{x}_i}^{{\tiny\between{1}{H+1}}}, w_1, u_1)$};
    \node (outdots1) at (8, -2) {$\cdots$};
    \node[xnode, fill=couleur3!30] (outLHpL) at (8,-3) {$\mbf{\hat{x}_i}^{{\tiny\between{L}{H+L}}} = m_\theta(\mbf{\hat{x}_i}^{{\tiny\between{L\shortminus1}{H+L\shortminus1}}}, u_L, w_L)$};

    \node (inv11) at (2.4, 0) {};
    \node (inv12) at (3.6, 0) {};
    \node (inv13) at (3.6, -0.5) {};
    \node (inv14) at (2.4, -0.5) {};
    \node (inv21) at (2.4, -1) {};
    \node (inv22) at (3.6, -1) {};
    \node (inv23) at (3.6, -1.5) {};
    \node (inv24) at (2.4, -1.5) {};
    \node (inv33) at (3.6, -2.5) {};
    \node (inv34) at (2.4, -2.5) {};
    \node (inv31) at (2.4, -3) {};
    \node (inv32) at (3.6, -3) {};
    
    \draw[arrow] (x0H.east) -- (inv11.west);
    \draw[arrow] (x1Hp1.east) -- (inv21.west);
    \draw[arrow] (xLHpL.east) -- (inv31.west);
    \draw[arrow] (inv12.east) -- (out0H.west);
    \draw[dashed] (out0H.south) |- (inv13.east);
    \draw[arrow, dashed] (inv14.west) -| (x1Hp1.north);
    \draw[dashed] (out1Hp1.south) |- (inv23.east);
    \draw[arrow, dashed] (inv24.west) -| (dots1.north);
    \draw[dashed] (outdots1.south) |- (inv33.east);
    \draw[arrow, dashed] (inv34.west) -| (xLHpL.north);
    \draw[arrow] (inv22.east) -- (out1Hp1.west);
    \draw[arrow] (inv32.east) -- (outLHpL.west);

    \draw[decorate, decoration={brace, amplitude=6pt}, thick]
      ($(out0H.east)+(1.1,0.3)$) --
      ($(outLHpL.east)+(0.3,-0.3)$)
      node[midway, right=8pt, align=center] {
          Compute\\[4pt]
          $\mathcal{L}_i\Bigl(\mbf{\hat{x}_i}^{{\tiny\between{H+1}{H+L}}}, \mbf{x_i}^{{\tiny\between{H+1}{H+L}}}\Bigr)$
        };
    \end{tikzpicture}}
    \caption{Training process using reference sequences $\mbf{x_i}$ and recursive calls to $m_\theta$ (Eq.~\ref{eq:training_loss}).}
    \label{fig:training}
\end{figure*}

\section{NMPC for Full-Scale Systems and Limitations of Learning-based Strategies}\label{sec:pb_statement}

\subsection{Background on MPC}
The principle of MPC is to use a dynamical model of the plant to predict and optimize its future behavior~\cite{schwenzer2021review}.
Most dynamical models can be represented as a set of continuous-time, first-order differential equations
\begin{equation}
    \dot{x}(t) = f\bigl( x(t), u(t), w(t) \bigr), \quad x(t_0) = x_0 \, ,
    \label{eq:ode}
\end{equation}
where $x(t) \in \R^{n_x}$ are the state variables, $u(t) \in \R^{n_u}$ are the control inputs, and $x_0 \in \R^{n_x}$ are the initial states. 
In this paper, the exogenous inputs $w(t) \in \R^{n_w}$ are assumed to be predictable.
Moreover, the model is assumed to include additional algebraic equations that stem from steady-state conditions, conservation laws, or geometric constraints
\begin{equation}
    \left\{
    \begin{aligned}
        \dot{x}_d(t) &= f\bigl( x_d(t), x_a(t), u(t), w(t) \bigr), \ \ x_d(t_0) = x_{d0} \,, \\
        0 &= g\bigl( x_d(t), x_a(t), u(t), w(t) \bigr) \,,
    \end{aligned}
    \right.
    \label{eq:dae}
\end{equation}
where $x_d(t) \in \R^{n_d}$ are the state variables and $x_a(t) \in \R^{n_a}$ are the algebraic variables.
Models based on differential-algebraic equations (DAEs) can offer a more natural representation of the plant than models based solely on ordinary differential equations (ODEs). 
However, DAE models are generally more challenging to simulate than ODE models, as they require consistent initialization of both differential and algebraic variables, tend to exhibit stiff behavior, and may have a high differential index~\cite{sjoberg2008optimal}. 
In addition, many practical applications involve switching logics or inequality constraints (e.g., dead zones or hysteresis), which further increase the computational burden. 
For these reasons, standard off-the-shelf DAE solvers often fail to meet industrial needs, motivating the use of dedicated in-house simulation algorithms. 

At each time step, the MPC actions are determined by repeatedly solving a finite-horizon OCP in real time. 
This OCP is typically formulated as
\begin{subequations} \label{eq:ocp}
    \begin{align}
        \min_{\ol{x}_d(\cdot),\ol{x}_a(\cdot),\ol{u}(\cdot)} \ \int\limits_{0}^{T} \ell &\bigl( \ol{x}_d(\tau),\ol{x}_a(\tau),\ol{u}(\tau), \ol{w}(\tau) \bigr) \odif{\tau} \,, \label{eq:ocp cost} \\
        \text{s.t. } \forall \tau \in [0,T], \ \dot{\ol{x}}_d &= f\bigl(\ol{x}_d(\tau),\ol{x}_a(\tau),\ol{u}(\tau), \ol{w}(\tau) \bigr) \,,\label{eq:ocp ode} \\
        0 &= g\bigl(\ol{x}_d(\tau),\ol{x}_a(\tau),\ol{u}(\tau), \ol{w}(\tau) \bigr) \,,\label{eq:ocp alg} \\
        x_{d0} &= \ol{x}_d(0) \,, \label{eq:ocp init} \\
        0 &\geq h\bigl(\ol{x}_d(\tau),\ol{x}_a(\tau),\ol{u}(\tau), \ol{w}(\tau) \bigr) \, , \label{eq:ocp ineq}
    \end{align}
\end{subequations}
where $T > 0$ is the length of the prediction horizon and $\ell$ is the cost function.
The equality constraints~\eqref{eq:ocp ode} and~\eqref{eq:ocp alg} enforce consistency with the model equations, while~\eqref{eq:ocp init} ensures that the initial model states match the current plant states $x_{d0}$. 
Additional constraints~\eqref{eq:ocp ineq}, whether equalities or inequalities, can be defined through the function $h$. 

In practice, the infinite-dimensional OCP~\eqref{eq:ocp} is transcribed into a finite-dimensional optimization problem using either single shooting, multiple shooting, or direct collocation~\cite{rao2009survey,diehl2009efficient}. 
However, this paper focuses on full-scale, industry-driven, long-horizon problems, where multiple shooting and direct collocation typically result in large optimization problems. 
Efficiently handling this increased complexity requires advanced implementation techniques, such as algorithmic differentiation and sparsity exploitation, which are often not supported by in-house DAE solvers.
For this reason, the OCP~\eqref{eq:ocp} is transcribed using single shooting, leading to the following Nonlinear Programming (NLP) formulation 
\begin{subequations} \label{eq:SS}
    \begin{align}
        \min_{\mbf{u}} \ &J(\mbf{u}|\Phi) \coloneqq \sum_{k=0}^{N-1} \ell\bigl(\Phi\bigl(t_{k+1}|x_{d0}, \mbf{u}, \mbf{w}\bigr), u_k\bigr) \,, \label{eq:cost_SS} \\
        \text{s.t. } \forall k &\in \llbracket 0,N-1 \rrbracket, \ 0 \geq h\bigl(\Phi\bigl(t_{k+1}|x_{d0},\mbf{u}, \mbf{w}\bigr), u_k\bigr)  \, .
    \label{eq:nlp}
    \end{align}
\end{subequations}
The model equations~\eqref{eq:ocp ode} -- \eqref{eq:ocp init} are solved externally by the in-house DAE solver $\Phi$, leaving only the control input sequence~$\mbf{u}$ as decision variable.
The optimal solution of the NLP~\eqref{eq:SS} is denoted by~$\mbf{u}^*$.
In the sequel, the state and algebraic variables are grouped into a single vector of model variables $x(t) \coloneqq (x_d(t), x_a(t))$ to simplify notation.
The symbol $\mbf{s}^{\between{a}{b}} \coloneqq \{s_a, \dots, s_b\}$ denotes the subsequence of $\mbf{s}$ corresponding to indices $0 \leq a \leq b \leq N$.

\subsection{Limitations of existing strategies}

\paragraph{Learning the controller} 
Direct learning of policies through Reinforcement Learning (RL) or Imitation Learning (IL) may appear promising to infer $\mbf{u}^*$ using $(x_{d0}, \mbf{w})$.
However, RL methods are computationally expensive~\cite{kaelbling1996reinforcement} and therefore unsuitable for simulation-intensive settings.
By contrast, IL generally enables faster and more stable training by mimicking the behavior of expert controllers designed for well-defined tasks.
Behavior cloning (BC), a simple supervised learning framework for implementing IL~\cite{foster2024behavior}, remains limited as it requires collecting large amounts of optimal control sequences $\mbf{u}^*$ for a large variety of $(x_{d0}, \mbf{w})$: if $\Phi$ is computationally expensive, generating such a dataset quickly becomes impractical.
More sophisticated methods based on IL, such as the one presented in~\cite{vaupel2020accelerating}, also highlight the heavy cost of generating the offline dataset.
In addition, methodologies like BC are prone to error amplification during long roll-outs as recursive predictions gradually diverge from the expert distribution, producing inaccurate trajectories in regions where the model cannot imitate correctly~\cite{foster2024behavior}.
Hence, in the presence of abrupt changes in $\mbf{w}$, $\ell$, or $h$, both RL and IL approaches may generate initial guesses that are far from optimal, leading to poor convergence under tight time budget.

\paragraph{Learning the dynamics} 
From another perspective, learning a surrogate model $m_\theta \approx \Phi$ has attracted increasing attention~\cite{chakraborty2021machine}, particularly for optimal control applications~\cite{de2021nonlinear}.
In this context, surrogate-based optimization replaces the expensive function $J(\mbf{u}) \coloneqq J(\mbf{u}|\Phi)$ with a lightweight approximation $J(\mbf{u}|m_\theta)$~\cite{han2012surrogate}.
Recent software tools such as \textit{L4CasADi}~\cite{salzmann2024learning} further facilitate the integration of such learned models directly within OCP solvers.
In the context of NMPC, embedding learned dynamics has emerged as a promising research direction~\cite{salzmann2023real}. 

\vspace{0.4em}%
\noindent In summary, existing warm-start strategies face a fundamental tension: approaches that directly learn $\mbf{u}^*$ suffer from prohibitive data generation costs when $\Phi$ is expensive, while the ones substituting $\Phi$ with a learned surrogate $m_\theta$ sacrifice the safety guarantees of the full-scale solver.
Learning the dynamics rather than the optimal policy -- and using $m_\theta$ solely to produce a warm-start that is then refined in the original problem~\eqref{eq:SS} -- offers a principled way to resolve this tension.
To the best of our knowledge, such a combination of surrogate dynamics learning and auxiliary optimization for structured warm-start generation in long-horizon NMPC has received limited attention, especially under last-minute changes in $\mbf{w}$, $\ell$, or $h$.

\section{Surrogate Dynamics Optimization for NMPC}\label{sec:methodology}

SDO aims to determine an initial guess for the NLP~\eqref{eq:SS} that accelerates NMPC convergence by exploiting a surrogate model of the dynamics. 
The latter acts as a lightweight substitute for the full-scale DAE solver $\Phi$, enabling the solution of a surrogate OCP with inexpensive dynamics prior to initializing the full NMPC problem.

\paragraph{Training phase (offline)} Let $m_\theta$ be the surrogate model parametrized by the vector $\theta$ (e.g., an arbitrary neural network). 
Assume access to a training set of $n_{\mbf{x}}$ reference trajectories from the full-scale DAE model. 
Given a control input sequence $\mbf{u_i}$ and an exogenous input sequence $\mbf{w_i}$, the sequence of model variables are collected as, for $i = 1, \dots, n_{\mbf{x}}$ and $x_i(t_0) \coloneqq x_{i0}$,
\begin{equation}
    \mbf{x_i}\!\coloneqq\!\bigl\{ x_{i0}, \Phi\bigl(t_{1}| x_{i0},\mbf{u_i}, \mbf{w_i}\bigr), \dots, \Phi\bigl(t_{N_i}| x_{i0}, \mbf{u_i}, \mbf{w_i}\bigr)\bigr\}.
\end{equation}
Let $H \geq 0$ be the context window. 
For a trajectory $i \in \llbracket 1, n_{\mbf{x}} \rrbracket$ and at time $k \in \llbracket H, N_i-1 \rrbracket$, the inputs $\mcal{I}_{i, k}$ and outputs $\mcal{O}_{i,k}$ of the surrogate model $m_\theta$ are defined by
\begin{equation}
    \mcal{I}_{i, k} \coloneqq \bigl(\mbf{x_i}^{\between{k\!\!-\!\!H}{k}}, (\mbf{u_i})_k, (\mbf{w_i})_k\bigr), \quad
    \mcal{O}_{i, k} \coloneqq \bigl(\mbf{x_i}\bigr)_{k+1}.
\end{equation}
To improve temporal coherence, the model is trained to predict sequences of length $L > 0$, allowing it to capture long-term dependencies. 
We denote by $m_\theta^{\between{1}{L}}(\mbf{x}^{[k-H, k]}, \mbf{u}, \mbf{w})$ the recursive prediction of the surrogate model over $L$ steps, using extra inputs $u_0, \dots, u_L$ and $w_0, \dots, w_L$ at appropriate times (see Fig.~\ref{fig:training}).

The distance between two sequences $\mbf{s_1}, \mbf{s_2} \in \R^L$ can be measured by a convex metric $d(\mbf{s_1}, \mbf{s_2})$ (e.g., mean squared error).
For each full trajectory $i \in \llbracket 1, n_{\mbf{x}} \rrbracket$, denote by
\begin{equation*}
\mbf{\hat{x}_i}^{\between{k\!+\!1}{k\!+\!L}}\coloneqq m_\theta^{[1:L]}\bigl(\mbf{x_i}^{\between{k\!\!-\!\!H}{k}}, \mbf{u_i}^{\between{k}{k\!+\!L\!-1}}, \mbf{w_i}^{\between{k}{k\!+\!L\!-1}}\bigr)
\end{equation*}
the predicted sequence and define the loss as
\begin{equation*}
    \mcal{L}_{i}(m_\theta) \coloneqq \sum_{k=H}^{N_i - L} d\bigl(\mbf{\hat{x}_i}^{\between{k\!+\!1}{k\!+\!L}}, \mbf{x_i}^{\between{k\!+\!1}{k\!+\!L}}\bigr) \, .
\end{equation*}
The loss over the whole training set is then given by
\begin{equation}\label{eq:training_loss}
    \mcal{L}(m_\theta) \coloneqq \frac{1}{n_{\mbf{x}}}\sum_{i=1}^{n_{\mbf{x}}} \mcal{L}_i(m_\theta) + \lambda \mcal{R}(m_\theta) \, ,
\end{equation}
where $\mcal{R}(m_\theta)$ is a regularization term that either penalizes model complexity or enforces physically feasible trajectories through a physics-informed prior~\cite{raissi2019physics, karniadakis2021physics}. 
The loss~\eqref{eq:training_loss} can be minimized using standard ML optimization algorithms, such as Adam~\cite{kingma2014adam}.

\paragraph{Inference phase: single shooting with the surrogate model} The surrogate model should be integrated in the optimization pipeline using single shooting, as in~\eqref{eq:SS}. 
We suggest removing hard constraints~\eqref{eq:nlp} and including them as soft constraints directly into the cost function for stability of the optimization routine, which is not guaranteed to be convex due to $m_\theta$. 
If $H = 0$, replacing the legacy integrator $\Phi$ with the surrogate model $m_\theta$ yields the lightweight NLP formulation
\begin{equation} \label{eq:SS_surrogate}
    \min_{\mbf{u}} \ \widetilde{J}(\mbf{u})\coloneqq \sum_{k=0}^{N-1} \Bigl[\ell \bigl(\mbf{\hat{x}}^{\between{1}{k}}, u_k\bigr) + \nu \bigl\|h^+\bigr(\mbf{\hat{x}}^{\between{1}{k}}, u_k\bigr)\bigr\| \Bigr]\,,
\end{equation} 
where $h^+$ denotes the component-wise positive part of $h$ and $\nu > 0$ is a penalty weight.
In practice, the more precise $m_\theta$ is, the larger $\nu$ should be set, to mimic hard constraints.
If $H \geq 1$, past measurements $\mbf{x}^{\between{\shortminus H}{0}} \coloneqq \{x(t_{\shortminus H}), \dots, x(t_0)\}$ should be stored and provided to $m_\theta$. 
The solution $\hat{\mbf{u}}$ of the lightweight NLP~\eqref{eq:SS_surrogate} should be obtained using a first-order optimization method (e.g., projected GD for box-shaped domains $\mcal{U}\subset \R^{n_u}$). 
Indeed, although the surrogate model mimics the DAE behavior, its evaluation may introduce approximation errors, resulting in a less smooth landscape. 
Thus, whereas second order optimization methods tend to take costly and larger steps using the Hessian matrix, gradient methods are more suited to noisy cost functions. 
The auxiliary optimization should be run for a small computational budget allocated before the main NLP resolution (e.g., $\leq 5\%$ of the total computational budget).

\paragraph{Warm-start phase: integration, evaluation and risk-management} Once the first optimization phase is completed, $\hat{\mbf{u}}$ is plugged into the original NLP~\eqref{eq:SS}. 
The quantitative quality of the warm-start can be assessed by comparing the cost function value or the constraints violation under a fixed time budget, or the number of iterations required to reach convergence. 
As this heuristic ends with an industrially trusted NMPC resolution, it remains inherently safe for real-life deployments in sensitive sectors. 
In the worst case, an inaccurate warm-start may only slow down convergence of the OCP. 

\paragraph{Bounding warm-start errors in Surrogate Dynamics Optimization} 
Variational analysis~\cite{rockafellar1998variational, bonnans2013perturbation} provides intuition that the optimal solution $\mbf{u}^*$ varies continuously with respect to the simulation scheme, under regularity assumptions. 
Informally, if $m_\theta$ approximates $\Phi$ uniformly well over the control space, $\hat{\mbf{u}}$ should lie in a neighborhood of $\mbf{u}^*$. 
Thus, the quality of the warm-start should be controlled by the approximation error. 
As an example, this conservative upper bound can be derived. 

\begin{proposition}\label{thm:upper_bound}
    Let $\mcal{X}_0$ and $\U$ be compact sets over 
    the state and control inputs spaces. For any $x_0\in\mathcal{X}_0$, let $\Phi^{\between{1}{N}}(x_0, \mbf{u}, \mbf{w}) \coloneqq \bigl\{\Phi\bigl(t_{1}|x_0, \mbf{u}, \mbf{w}\bigr), \dots, \Phi\bigl(t_{N}|x_0, \mbf{u}, \mbf{w}\bigr)\bigr\}$.
    Assume there exists $M > 0$ such that
    \begin{equation*}\label{eq:thm_eq_1}
        \max_{\mbf{u} \in \U} \|\Phi^{\between{1}{N}}(x_0, \mbf{u}, \mbf{w}) - m_\theta^{\between{1}{N}}(x_0, \mbf{u}, \mbf{w})\| \leq M \, .
    \end{equation*}
    Consider the two OCPs
    \begin{equation*}
        \begin{aligned}
            \mbf{u}^* \coloneqq \argmin_{\mbf{u}\in \U} J(\mbf{u} |\Phi), \quad
            \hat{\mbf{u}} \coloneqq \argmin_{\mbf{u}\in \U} J(\mbf{u} |m_\theta) \, .
        \end{aligned}
    \end{equation*}
    Assume that $\forall \mbf{u}\in\mcal{U}, \, \mbf{x} \mapsto J_\mbf{x}(\mbf{u}) \coloneqq \sum_{k=0}^{N-1} \ell(\mbf{x}_k, u_k)$ is $K_{\!J}$-Lipschitz, and that
    $J(\cdot|\Phi)$ is differentiable and $\mu$-strongly convex. 
    Then we have the following bound on $\|\mbf{u}^* - \hat{\mbf{u}}\|$:
    \begin{equation*}\label{eq:theoretical_bound}
        \|\mbf{u}^* - \hat{\mbf{u}}\| \leq 2\sqrt{\frac{K_{\!J} M}{\mu}} \, .
    \end{equation*}
\end{proposition}

\begin{remark}
    Though $J(\cdot|\Phi)$ is unlikely to be strongly convex in general, Prop.~\ref{thm:upper_bound} describes the central mechanism justifying our approach: increasing the precision of the surrogate model directly impacts the quality of the initialization.
    Without the assumption of strong convexity of $J(\cdot|\Phi)$, one can obtain the bound $J(\hat{\mbf{u}}|\Phi) - J(\mbf{u}^*|\Phi) \leq 2K_J M$ on the objective value gap, which can be useful in non-convex settings.
\end{remark}

This property, while established under simplified and conservative assumptions, is empirically supported by Fig.~\ref{fig:results_2} which shows a monotone relationship between surrogate accuracy and warm-start quality.
However, in complex settings, the lack of theoretical guarantees prevents us from removing the second optimization process, as it acts like a safety filter recovering feasibility with respect to the original constraints and model~\cite{hsu2023safety}.

\section{Industrial case study}\label{sec:experiments}

The following case study is deliberately set in a full-scale industrial context, a PWR under load-following control, rather than on a toy dynamical system. 
This choice reflects the core motivation of this work: SDO targets settings where the computational bottleneck is the simulation of a high-dimensional, stiff DAE model, a regime in which standard benchmarks would fail to highlight the practical limitations of existing warm-start strategies.

\subsection{System modeling of the PWR}

To address the challenges posed by the evolving energy mix, recent works have proposed NMPC-based systems to enhance the load-following capabilities of PWRs~\cite{dupre2021enhanced}. 
Monitoring a PWR is a complex task, and the OAPS solution~\cite{dupre2025oaps} provides control room operators with real-time recommendations derived from a set of optimal strategies. 
However, accurately modeling a PWR for long-term planning strategies, with prediction horizons extending up to 24~hours, is a major challenge due to the required processing effort.
Reducing the computation time of these demanding strategies is essential to ensure that the tool remains responsive to last-minute grid requests.
In such cases, standard shift-initialization techniques can slow down the convergence of the OCP, since the newly requested scenario may differ completely from the original one. 
The considered one-dimensional PWR model, used to predict the long-term strategies of the OAPS solution, is described by the following set of DAEs
\begin{subequations} \label{eq:PWR model}
        \begin{align}
             \odv{\I(t)}{t} &= \gamma_{\sss\I} \Pn(t) - \lambda_{\sss\I} \I(t) \,, \\ 
            \odv{\Xe(t)}{t} &= \gamma_{\sss\Xe} \Pn(t) + \lambda_{\sss\I} \I(t) - \bigl(\lambda_{\sss \Xe} + \sigma_{\sss\Xe} \Pn(t) \bigr) \!\odot\! \Xe(t) \,, \\
            \odv{\h(t)}{t} &= u(t) \,, \\
            0 &= \rho(t) \odot \Pn(t) + D M_{\mrm{exch}} \Pn(t) \,, \\
            0 &= \sum_{i=1}^{\nz} \Pn_i(t) - w(t) \,, \\
            0 &= \Tin(t) - \Tref \bigl( w(t) \bigr) \,, 
        \end{align}
\end{subequations}
\noindent where $n_z \in 2\N^*$ is the even number of axial nodes in the reactor core, $\odot$ represents the element-wise product, and $M_{\mrm{exch}} \coloneqq \mrm{diag}_{\nz}\!(-1, -2, \dots, -2, -1) + I_{\nz}^{+} + I_{\nz}^{-} \in \R^{n_z \times n_z}$ with $I_{\nz}^+$ and $I_{\nz}^-$ denoting the unit super-diagonal and sub-diagonal matrices, respectively.
Each node is characterized by its iodine concentration $\I_j$, xenon concentration $\Xe_j$, and fission power $\Pn_j$.
The inlet reactor core temperature $\Tin$ follows the reference temperature program $\Tref$, which depends on the turbine power.
The turbine power is treated as a predictable exogenous input $w(t) \in \R$, since the electrical load schedule is communicated in advance by the transmission system operator.
Power balance between the turbine and the reactor is maintained by adjusting the boron concentration $\Cb$ in the primary coolant and the cumulative position $\h$ of the control rod groups.
These variables, along with the inlet temperature, xenon concentrations, and node power levels, affect the growth rate of the fission chain reaction through the nonlinear terms
\begin{equation*}
\rho_j(t)\coloneqq\rho_j\bigl(\Pn_1(t),\ldots,\Pn_j(t),\Tin(t),\Xe_j(t),\h(t),\Cb(t) \bigr)\, ,
\end{equation*}
for $j = 1, \ldots, \nz$, revealing the complex feedbacks and strong couplings that make reactor core control challenging~\cite{reuss2003precis, stacey2018nuclear}.
The state and algebraic variables of the PWR model~\eqref{eq:PWR model} are respectively given by
\begin{equation*}
    \begin{aligned}
         x_d(t) &= \bigl[ \I_1(t), \, \ldots, \, \I_{\nz}(t), \, \Xe_1(t), \ldots, \, \Xe_{\nz}(t), \, \h(t) \bigr]^\top, \\
         x_a(t) &= \bigl[ \Pn_1(t), \, \ldots, \, \Pn_{\nz}(t), \, \Cb(t), \, \Tin(t) \bigr]^\top \ .
    \end{aligned}
\end{equation*}
The control input $u(t) \in \R$ corresponds to the speed of the control rod groups.

\subsection{Experimental setup}\label{sec:exp_setup}

\begin{figure*}
    \centering
    \begin{subfigure}[t]{0.32\linewidth}
        \centering
        \includegraphics[width=\linewidth]{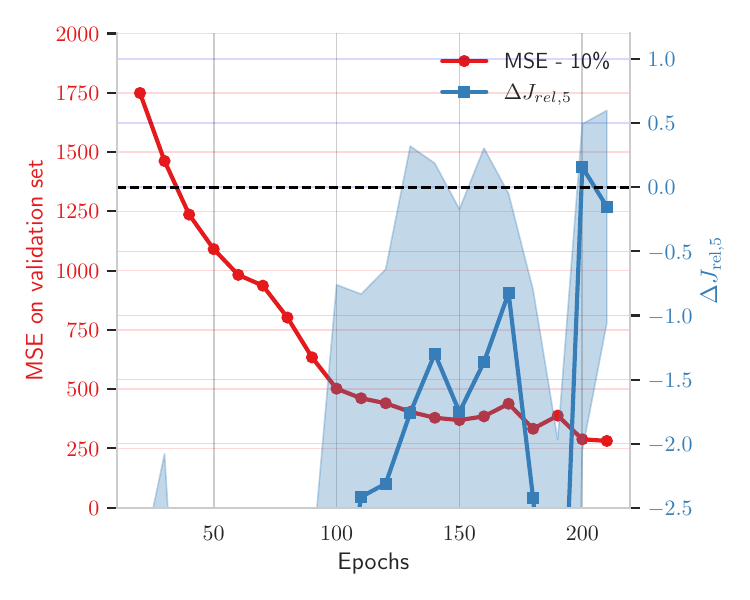}
        \caption{Training on 10\% of the dataset.
        The surrogate fails to provide reliable acceleration:
        the median $\Delta J_{\mathrm{rel},5}$ remains close to 0
        or negative throughout training, indicating insufficient
        data to learn dynamics useful for warm-starting.}
        \label{fig:results_2_10}
    \end{subfigure}
    \hfill
    \begin{subfigure}[t]{0.32\linewidth}
        \centering
        \includegraphics[width=\linewidth]{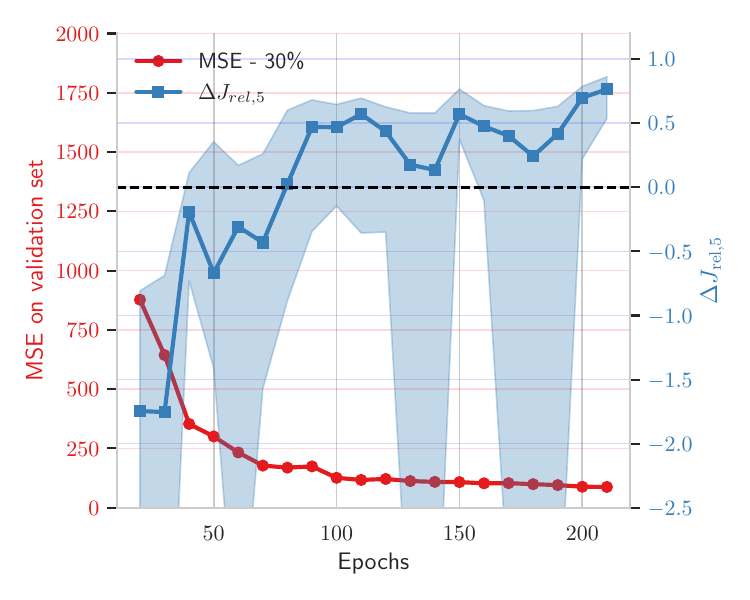}
        \caption{Training on 30\% of the dataset.
        Significant acceleration is achieved at convergence,
        demonstrating that the framework does not require
        exhaustive data collection. The MSE--$\Delta J_{\mathrm{rel},5}$
        correlation emerges, though with higher variance than
        with the full dataset.}
        \label{fig:results_2_30}
    \end{subfigure}
    \hfill
    \begin{subfigure}[t]{0.32\linewidth}
        \centering
        \includegraphics[width=\linewidth]{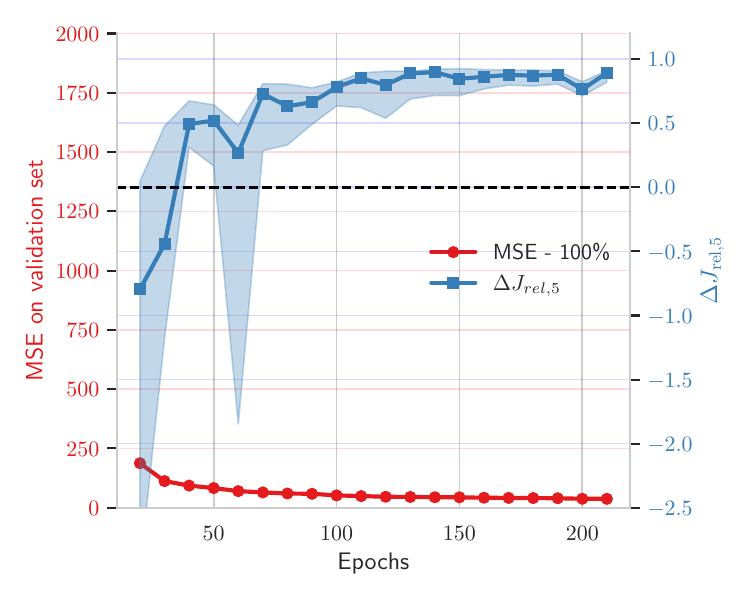}
        \caption{Training on the full dataset (10000 time-series).
        A clear monotone relationship between validation MSE
        and $\Delta J_{\mathrm{rel},5}$ emerges after $\sim$50 epochs,
        confirming that surrogate accuracy directly drives
        warm-start quality. The stopping criterion should
        account for the relative gain, not MSE alone.}
        \label{fig:results_2_full}
    \end{subfigure}
    \caption{Comparison between surrogate accuracy (MSE on validation set)
    and the relative gain $\Delta J_{\mathrm{rel},5}$ (defined in Eq.~\eqref{eq:acceleration_metric}) under abrupt changes in $\mbf{w}^{\between{k}{k+L}}$ across three training set sizes.
    The solid line shows the median $\Delta J_{\mathrm{rel},5}$ and the shaded area indicates interquartile range (Q25--Q75).
    Surrogate optimization accounts for up to $5\%$ of the total computational budget in all experiments.}
    \label{fig:results_2}
\end{figure*}

\paragraph{Validation scenario} Regarding the OCP, the prediction horizon is set to $T=24$~h, with a simulation time step of $T_s = 10$~min. 
The NMPC problem is solved in closed-loop at each time step.
The piecewise-constant control trajectory $u(t)$ is therefore parametrized by $N=L=144$ decision variables. 
The considered optimization strategies are defined by the instantaneous cost
\begin{subequations} \label{eq:strategies}
    \begin{align}
        \ell_{\sss\Pn}(\mbf{x}, \mbf{u}) &= \bigl(\AO(\Pn(t_{k+1})) - \AO_{\sss\Pn}^{\mrm{ref}}\bigr)^2, \label{eq:target_P} \\ 
        \ell_\mrm{b}(\mbf{x}, \mbf{u}) &= \bigl(\Cb(t_{k+1}) - \Cb(t_{k})\bigr)^2, \label{eq:target_bore}
    \end{align} 
\end{subequations}
\noindent where $\AO(\Pn) \coloneqq (\sum_{j=1}^{\sfrac{\nz}{2}} \Pn_j - \sum_{j=\sfrac{\nz}{2}+1}^{\nz} \Pn_j)/(\sum_{j=1}^{\nz}\Pn_j)$ represents the Axial Offset, i.e., the axial imbalance of power $\Pn$ between the top and bottom halves of the reactor core. 
The cost function~\eqref{eq:target_P} mitigates axial power oscillations, while~\eqref{eq:target_bore} reduces effluent production.
The goal is to return a sequence of control rod group speeds $\mbf{u}$ minimizing the integral cost $J(\mbf{u})$ of instantaneous cost $\ell_{\sss\Pn}$ or $\ell_\mrm{b}$.
In both setups, nonlinear hard constraints of the form $\AO_\mathrm{min} \leq \AO(\Pn(t_{k+1})) \leq \AO_\mathrm{max}$ are mandatory at all times for safety reasons.

\begin{table}[t]
\centering
\small
\setlength{\tabcolsep}{5pt}
\begin{tabular}{lcc}
\toprule
\textbf{Architecture} & MSE & Median $\Delta J_\mathrm{rel, 5}$ \\
\cmidrule(lr){1-3}
MLP & 2038.9 & -0.71 \\
RNN & 1697.1 & -0.69 \\
Transformer & \textbf{34.8} & \textbf{0.85} \\
\bottomrule
\end{tabular}
 
\caption{Performance comparison of different architectures (optimized with \textit{Optuna} -- 100 trials, max 30000 parameters). Results are computed on the validation set for trajectory reconstruction (MSE) and the target metric $\Delta J_\mathrm{rel, 5}$ (see Sec.~\ref{sec:results}) under abrupt changes in $\mathbf{w}^{\between{k}{k+L}}$. Small MLP and RNN architectures fail to properly capture the system dynamics, yielding no acceleration.
\vspace{-5pt}}
\label{tab:benchmark_arch}
\end{table}

\paragraph{Designing the surrogate model} The surrogate model aims to learn the PWR dynamics modeled by $\nz = 6$ axial nodes, yielding 13 state variables and 8 algebraic variables. 
To achieve this, the surrogate model learned over $n_{\mbf{x}} = 10000$ transients of $24$~h, randomly divided in $80\%$ for training and $20\%$ for validation. 
To select the model, a comparison of several classic neural network architectures (multi-layer perceptron, recurrent neural networks, and transformer-based models) was conducted. 
Hyperparameters were optimized with \textit{Optuna}~\cite{akiba2019optuna} over 100 trials, with the total number of model parameters constrained to be at most 30000.
The presented results are obtained with a simple transformer-based encoder-decoder architecture, which managed better results (see Tab.~\ref{tab:benchmark_arch}) as it tamed better temporal dependencies and showed more stability on long-horizon roll-outs.
Small MLP and RNN architectures, despite hyperparameter optimization, fail to capture the long-horizon temporal dependencies of the PWR dynamics, producing surrogate landscapes that mislead the auxiliary optimization. 
This underlines that surrogate quality is a prerequisite for the framework to be beneficial, consistent with Fig.~\ref{fig:results_2}.
The surrogate model is trained with a context window of $H = 1$ (i.e., measures from $t = -T_s$ and $t=0$ mins) until convergence over the prediction MSE.
This minimal context was found sufficient for the Transformer architecture to capture the memory of the PWR dynamics. 
Larger values of $H$ were tested but yielded no significant improvement.
Overall, the surrogate model is 100-500$\times$ faster to evaluate than the in-house DAE solver. 
All experiments were conducted on a NVIDIA Tesla V100-SXM2-32GB GPU using \textit{PyTorch}.

\subsection{Results}\label{sec:results}

\definecolor{darkgreen}{RGB}{0,0,0}
\definecolor{darkred}{RGB}{0,0,0}
\begin{table*}[t]
\centering
\small
\setlength{\tabcolsep}{4.5pt}
 
\begin{tabular}{llcccccccc}
\toprule
\multirow{4}{*}{\textbf{Statistic}} &
\multirow{4}{*}{\textbf{Method}} &
\multicolumn{8}{c}{\textbf{Relative improvement } $\Delta J_{\mathrm{rel},n}$} \\
 
\cmidrule(lr){3-10}
& &
\multicolumn{2}{c}{$n=5$} &
\multicolumn{2}{c}{$n=10$} &
\multicolumn{2}{c}{$n=15$} &
\multicolumn{2}{c}{$n=20$} \\
 
\cmidrule(lr){3-4}
\cmidrule(lr){5-6}
\cmidrule(lr){7-8}
\cmidrule(lr){9-10}
& & $\mbf{w}^{\between{k}{k+L}}$ & $J(\mbf{u})$
& $\mbf{w}^{\between{k}{k+L}}$ & $J(\mbf{u})$
& $\mbf{w}^{\between{k}{k+L}}$ & $J(\mbf{u})$
& $\mbf{w}^{\between{k}{k+L}}$ & $J(\mbf{u})$ \\
\midrule
 
\multirow{4}{*}{Worst}
& Shift Init
& \textcolor{darkred}{-2.65} & \textcolor{darkred}{-1.84}
& \textcolor{darkred}{-2.57} & \textcolor{darkred}{-12.05}
& \textcolor{darkred}{-3.48} & \textcolor{darkred}{-21.20}
& \textcolor{darkred}{-3.55} & \textcolor{darkred}{-22.15} \\
& BC
& -2.02 & -0.23
& -2.11 & -0.45
& -2.83 & -0.42
& -8.31 & -0.38 \\
& SDO (1\%)
& -0.34 & 0.29
& -0.31 & 0.11
& -0.01 & 0.00
& -5.10 & -0.06 \\
& SDO (5\%)
& \textbf{-0.15} & \textbf{0.30}
& \textbf{0.35} & \textbf{0.13}
& \textbf{0.35} & \textbf{0.04}
& \textbf{-0.70} & \textbf{-0.01} \\
\midrule
 
\multirow{4}{*}{Q25}
& Shift Init
& 0.00 & 0.06
& \textcolor{darkred}{-0.08} & \textcolor{darkred}{-0.80}
& \textcolor{darkred}{-0.29} & \textcolor{darkred}{-2.12}
& \textcolor{darkred}{-0.10} & \textcolor{darkred}{-4.35} \\
& BC
& 0.55 & 0.09
& 0.30 & 0.13
& 0.15 & 0.11
& 0.20 & 0.15 \\
& SDO (1\%)
& 0.81 & 0.37
& 0.72 & 0.33
& 0.70 & 0.30
& 0.66 & 0.29 \\
& SDO (5\%)
& \textbf{0.83} & \textbf{\textcolor{darkgreen}{0.43}} 
& \textbf{0.77} & \textbf{\textcolor{darkgreen}{0.41}} 
& \textbf{0.77} & \textbf{\textcolor{darkgreen}{0.38}} 
& \textbf{0.74} & \textbf{\textcolor{darkgreen}{0.34}} \\
\midrule
 
\multirow{4}{*}{Median}
& Shift Init
& 0.31 & 0.29
& 0.22 & 0.20
& 0.20 & 0.18
& 0.21 & 0.19 \\
& BC
& 0.72 & 0.36
& 0.64 & 0.34
& 0.54 & 0.31
& 0.50 & 0.33 \\
& SDO (1\%)
& 0.86 & 0.57 
& 0.87 & 0.47 
& \textbf{0.88} & 0.45 
& 0.85 & 0.40\\
& SDO (5\%)
& \textbf{\textcolor{darkgreen}{0.91}} & \textbf{\textcolor{darkgreen}{0.61}} 
& \textbf{\textcolor{darkgreen}{0.88}} & \textbf{\textcolor{darkgreen}{0.55}} 
& \textbf{\textcolor{darkgreen}{0.88}} & \textbf{\textcolor{darkgreen}{0.52}} 
& \textbf{\textcolor{darkgreen}{0.86}} & \textbf{\textcolor{darkgreen}{0.48}} \\
\midrule
 
\multirow{4}{*}{Q75}
& Shift Init
& \textcolor{darkgreen}{0.70} & 0.36
& \textcolor{darkgreen}{0.51} & 0.31
& \textcolor{darkgreen}{0.55} & 0.34
& \textcolor{darkgreen}{0.53} & 0.35 \\
& BC
& 0.80 & 0.83
& 0.72 & 0.78
& 0.65 & 0.73
& 0.60 & 0.68 \\
& SDO (1\%)
& \textbf{0.91} & 0.93 
& \textbf{0.90} & 0.93
& 0.91 & 0.92 
& 0.91 & 0.92 \\
& SDO (5\%)
& \textbf{\textcolor{darkgreen}{0.91}} & \textbf{\textcolor{darkgreen}{0.94}} 
& \textbf{\textcolor{darkgreen}{0.90}} & \textbf{\textcolor{darkgreen}{0.94}} 
& \textbf{\textcolor{darkgreen}{0.93}} & \textbf{\textcolor{darkgreen}{0.94}} 
& \textbf{\textcolor{darkgreen}{0.92}} & \textbf{\textcolor{darkgreen}{0.93}} \\
\bottomrule
\end{tabular}

\caption{
Statistics of the relative gain $\Delta J_{\mathrm{rel},n}$ defined in~\eqref{eq:acceleration_metric} after $n=5,10,15,20$ gradient steps under abrupt changes in either $\mbf{w}^{\between{k}{k+L}}$ or $J(\mbf{u})$ over 20~runs.
SDO ($k$\%) indicates the fraction of computational time budget allocated to surrogate optimization.
All methods were given a budget of 20 full-scale iterations ($\approx 20$ min of execution time).
Values closer to 1 indicate effective warm-starting, whereas negative values indicate performance deterioration.
\textbf{Bold values} highlight the best method.
}
 
\label{tab:benchmark_costs}
\end{table*}

\begin{table}[t]
\centering
\small
\setlength{\tabcolsep}{5pt}
\begin{tabular}{lcccccccccc}
\toprule
\multirow{4}{*}{\textbf{Method}} &
\multicolumn{6}{c}{\textbf{Constraints violations at iteration $n$}} \\
 
\cmidrule(lr){2-7}
&
\multicolumn{2}{c}{$n=0$} &
\multicolumn{2}{c}{$n=5$} &
\multicolumn{2}{c}{$n=10$} \\
 
\cmidrule(lr){2-3}
\cmidrule(lr){4-5}
\cmidrule(lr){6-7}
& $\E_{h^+}$ & $\P_{h^+}$
& $\E_{h^+}$ & $\P_{h^+}$
& $\E_{h^+}$ & $\P_{h^+}$ \\
\midrule
 
Cold-start
& $2\cdot 10^{4}$ & 100\%
& 60.6 & 15\%
& \textbf{0} & \textbf{0\%} \\
SI ($\ne J(\mbf{u})$)
& \textbf{31.8} & 50\%
& \textbf{0} & \textbf{0\%}
& \textbf{0} & \textbf{0\%} \\
SI ($\ne \mbf{w}$)
& $5\cdot 10^3$ & 100\%
& 0.2 & 5\%
& \textbf{0} & \textbf{0\%}\\
BC
& $3.5\cdot 10^2$ & 50\%
& \textbf{0} & \textbf{0\%}
& \textbf{0} & \textbf{0\%} \\
SDO (1\%)
& $6.5\cdot 10^2$ & 45\%
& 1.5 & 5\%
& \textbf{0} & \textbf{0\%} \\
SDO (5\%)
& $6.7\cdot 10^2$ & \textbf{35\%}
& 1.6 & 5\%
& \textbf{0} & \textbf{0\%} \\

\bottomrule
\end{tabular}
 
\caption{
Average constraint violation $\E_{h^+} \coloneqq \E(\| h^+ \|)$, and probability of constraint violation $\P_{h^+}\coloneqq\P(\| h^+ \| > 0)$, for different warm-start strategies at iterations $n = 0, 5,10$ of the optimization process.
Shift-initialization with abrupt changes in $J(\mbf{u})$ performs best because the constraints are not modified. 
By contrast, SDO generally enforces feasibility faster when $\mbf{w}^{\between{k}{k+L}}$ changes.
\vspace{-5pt}}
\label{tab:benchmark_constraints}
\end{table}

\begin{figure*}
    \centering
    \includegraphics[width=\linewidth]{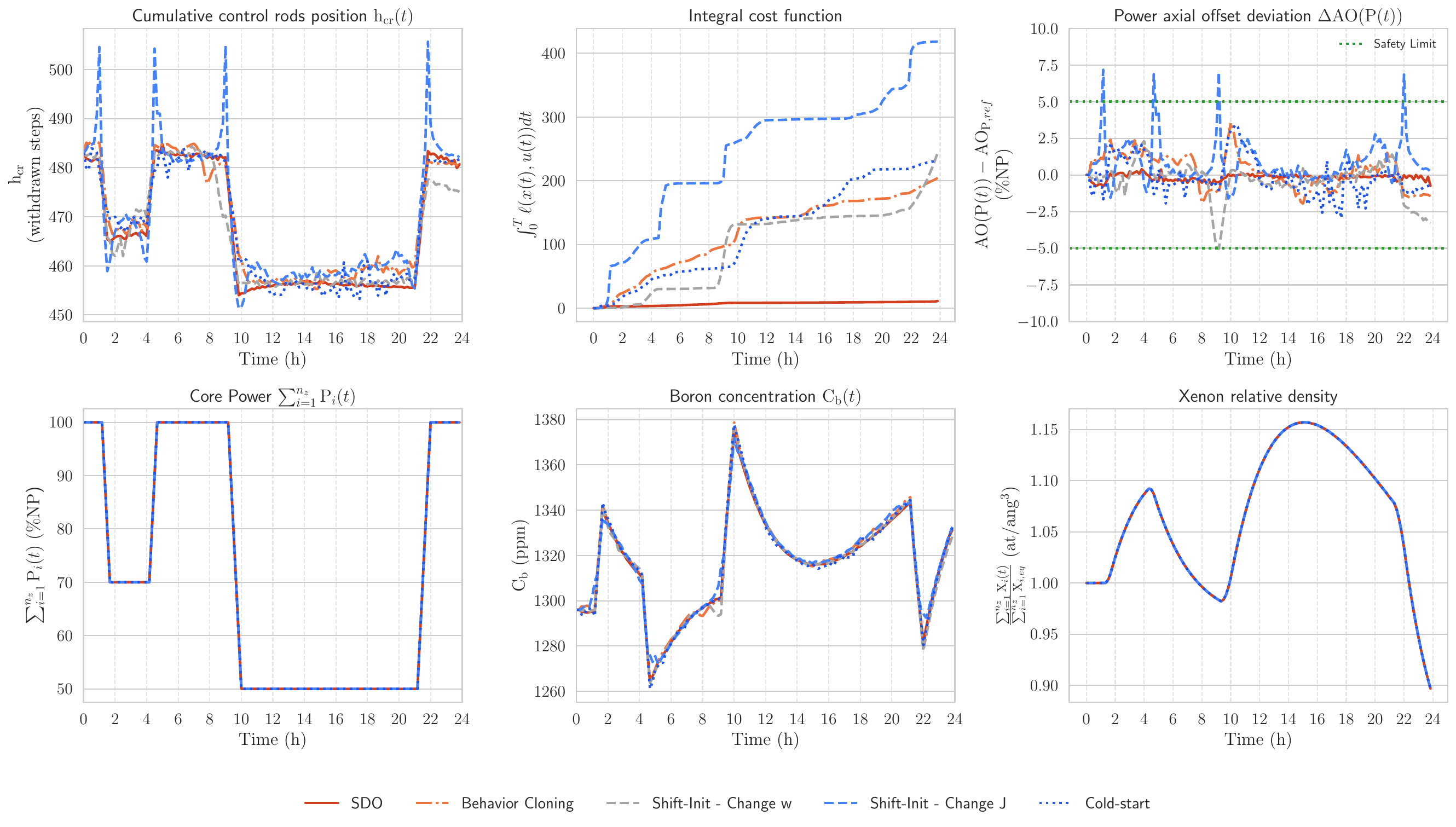}
    \caption{Comparison of relevant state and algebraic variables during a reference transient using the control input $\mbf{u}$ obtained after 10~NMPC iterations. Results are shown for cold-start (blue, dotted), shift initialization starting from a different $J(\mbf{u})$ (blue, dashed), shift initialization starting from a different $\mbf{w}^{\between{k}{k+L}}$ (grey, dashed), behavior cloning (orange, dashdotted) and SDO (solid, red). The optimized cost corresponds to the integral cost function over 24h (top-center graph).}
    \label{fig:example}
\end{figure*}

We compared the solutions obtained using 4 different warm-start strategies: cold-start, shift-initialization (SI), behavior cloning (BC), and SDO. 
All experiments were conducted in scenarios involving sudden changes, where either the future exogenous inputs $\mbf{w}^{\between{k}{k+L}}$ or the cost function $J(\mbf{u})$ is modified at the last minute (i.e., at time $k$).
The target is to react as quickly as possible in order to find a new satisfactory control input strategy.

\paragraph{Baselines}
Cold-start procedure consists in control inputs initialized at $\mbf{u}_{\text{cold}} \coloneqq \{ 0, \dots, 0 \}$), while shift-initialization procedure consists in control inputs initialized through an optimization over a different strategy $J(\mbf{u})$ or $\mbf{w}^{\between{k}{k+L}}$ at previous timestep.  
For BC, we optimized and trained an MLP using \textit{Optuna} to directly output the optimal control inputs knowing $(\mbf{x}^{\between{k-H}{k}}, J(\mbf{u}), \mbf{w}^{\between{k}{k+L}})$, with $H = 1$.
Other architectures such as RNN or Transformers of similar size were also tested but did not show any significant improvement as no recursive predictions are needed for this task.
However, as mentioned in Sec.~\ref{sec:Introduction}, the offline cost to train the BC network is very heavy as it requires a dataset of previously computed optimal control sequences $\mbf{u}^*$ for a variety of different $(\mbf{x}^{\between{k-H}{k}}, J(\mbf{u}), \mbf{w}^{\between{k}{k+L}})$. 
We want to highlight that the limit of a BC controller is not its ability to infer the mapping $(\mbf{x}^{\between{k-H}{k}}, J(\mbf{u}), \mbf{w}^{\between{k}{k+L}}) \mapsto \mbf{u}^*$, which achieved satisfactory accuracy on the training set, but the high computational cost needed to produce a large amount of optimal control inputs $\mbf{u}^*$.
To generate the training set for BC, 300 full-scale OCPs were solved. 
This represents about 1 million trajectories, which corresponds to $\approx$100 hours of cumulated computational time.
In comparison, the computational time spent generating the training set for the surrogate model was less than 1 hour.
Despite the large computational time spent generating the training set, the dataset of optimal trajectories proved insufficient for the BC network to correctly generalize.
\paragraph{Comparison methodology}
For fairness of the online comparison, we used a constant computational budget characterized by fixed amount of gradient descent iterations (due to the high computational cost of $J$, second order optimization methods are not viable). 
The metric used is defined by 
\begin{equation}\label{eq:acceleration_metric}
    \Delta J_\mathrm{rel, n}\!=\!\frac{J(\mbf{u}_n\!\mid\!\mbf{u}_0 = \mbf{u}_{\mathrm{cold}}) -
      J(\mbf{u}_n\!\mid\!\mbf{u}_0\!=\!\mbf{u}_{\mathrm{strat}})}
     {J(\mbf{u}_n\!\mid\!\mbf{u}_0\!=\!\mbf{u}_{\mathrm{cold}})},
\end{equation} which represents the relative improvement from a given method $\mathrm{strat}$ compared to the cold-start procedure at iteration $n$ of the optimization process.
The results obtained over a validation set of transients with different strategies indicate that SDO achieves on average higher performance than its competitors (quantitative results in Tab.~\ref{tab:benchmark_costs} and Tab.~\ref{tab:benchmark_constraints}, example in Fig.~\ref{fig:example}).
Moreover, the results demonstrate that a very small computational budget (e.g., 1\% of the total budget, $\approx$ 20 s) is required to achieve significant acceleration.
This improvement facilitates real-time optimization and a faster response to the electricity regulator. 
Fig.~\ref{fig:results_2} further investigates the sensitivity of the warm-start quality to both training duration and dataset size. 
With the full dataset (100\%), a clear monotone relationship between validation MSE and $\Delta J_\mathrm{rel, 5}$ emerges after approximately 50 epochs, confirming that surrogate accuracy directly drives warm-start quality. 
Crucially, using only 30\% of the training data (3000 trajectories) still yields a median $\Delta J_\mathrm{rel, 5}$ above 0 at convergence, demonstrating that SDO does not require exhaustive data collection. 
In contrast, 10\% of the data proves insufficient for reliable acceleration. 
This data efficiency is a key advantage over BC-based approaches, which require a comparably large dataset of optimal trajectories.



\section{Conclusion and future works}\label{sec:conclusion}

In this paper, we describe SDO, a warm-start strategy based on a ML surrogate model designed to accelerate full-scale NMPC problems. 
The surrogate model acts as an auxiliary simulation scheme integrated into an optimization routine. 
The output of this lightweight problem is used to initialize the original NLP. 
SDO is reproducible on any NMPC problem, and numerical results demonstrate a significant performance gain for a fixed time budget compared to standard warm-start techniques on systems that are costly to integrate. 

This approach suggests that surrogate models can be leveraged as NMPC accelerators rather than as substitutes. 
This supports the view that data-driven and model-based control methods are complementary, even for sensitive applications. 
In our specific application, this methodology could allow electricity regulators to optimize trajectories over multiple days, which is currently not feasible.
Future directions include how to adapt the presented procedure in other contexts including stochastic ones, as we believe our contribution is not limited to the nuclear industry. 
Moreover, a larger comparison of warm-start methods for NMPC problems under a variety of different contexts would be valuable to establish practical guidelines for industrial applications. 
Finally, extending the theoretical framework under less restrictive assumptions constitutes a promising direction.



\bibliographystyle{IEEEtran}
\bibliography{refs.bib}             

\end{document}